\documentclass[reqno]{amsart}
\usepackage{amsmath}%
\usepackage{amsfonts}%
\usepackage{amssymb}%
\usepackage{ifpdf}
\ifpdf
\usepackage{graphicx}
\else
\usepackage[draft]{graphicx}
\fi
\usepackage{geometry}
\usepackage{subfig}
\begin{document}
\title{A Method for Zeroing-In on $\text{Re}\,\zeta(\sigma+it)<0$ in the half-plane $\sigma>1$}
\author{Dominic C. Milioto}
\date{January 12, 2010}
\subjclass[2000]{Primary 33F05,30-02; Secondary 65K05, 68W40} %
\keywords{Zeta function, congruences, modular arithmetic}%

\begin{abstract}
This paper describes a search algorithm to locate values of $t$ where the real part of the Riemann zeta function, $\zeta(\sigma+it)$, is negative for $\sigma>1$.  The run-time to execute the search is much less than a brute-force approach and relies on certain symmetries of congruence equations related to the zeta function.
\end{abstract}
\maketitle

\section{Introduction}
\label{intro}
The goal of this study is three-fold:
\begin{enumerate}
\item Construct an algorithm to approximate a solution to a system of congruence equations pertaining to the zeta function.
\item Use the approximations in item (1) to locate regions along the line $\sigma=1$ where the real part of $\zeta(\sigma+it)$ may become negative.
\item Attempt to locate and verify instances where the real part of this function  becomes negative along this line by a method other than trial and error.
\end{enumerate}

The value of $\zeta(s)$ for $s=\sigma+it$ in the half-plane $\sigma>1$ is commonly given by it's Dirichlet series
\begin{equation}
\zeta(s)=\sum_{n=1}^{\infty}\frac{1}{n^s},\quad \text{Re}(s)>1
\end{equation}
or it's Euler Product form
\begin{equation}
\zeta(s)=\prod_{p}\frac{1}{1-1/p^s},\quad \text{Re}(s)>1
\label{eqn001}
\end{equation}
where $p$ is over the set of primes.
From these two expressions, little can be deduced about the range of $\zeta(s)$.  However, Titchmarsh in \cite{Titchmarsh2}, provides an analysis which gives an upper bound.  We first take the logarithm of both sides of \eqref{eqn001} to obtain

\begin{equation}
\log\zeta(\sigma+it)=-\sum_{n=1}^{\infty} \log\left(1-p_n^{-s}\right).
\label{eqn002}
\end{equation}

For $s=\sigma_0+it$ with $\sigma_0>1$, let:
$$U=\{u: u=\log\zeta(\sigma_0+it)=-\sum_{n=1}^{\infty} \log\left(1-p_n^{-\sigma_0} e^{-it\log(p_n)}\right)$$
and
$$V=\{v: v=\Phi(\sigma_0,\theta_1,\theta_2,\cdots)=-\sum_{n=1}^{\infty} \log\left(1-p_n^{-\sigma_0} e^{-i\theta_n\log(p_n)}\right)$$
for independent variables $(\theta_1,\theta_2,\cdots)$ with $0\leq \theta_n \leq 2\pi$.
Note first that $U\subset V$ since $t\log p_n \equiv \theta_n \text{mod}\; 2\pi$.

Using the set $V$, Titchmarsh placed approximate bounds on the range of the set $U$ as a function of $\sigma$ for $\log \zeta(s)$ in the half-plane $\sigma > 1$.  These ranges are in the form of circles with centers at $1/2\log \zeta(2 \sigma)$ given by \eqref{eqnveqn}.
\begin{equation}
U = \left\{u: u\in \text{D}\left(1/2\log \zeta(2\sigma), \frac {1}{2}\log\frac {\zeta^2(\sigma)}{\zeta(2\sigma)}\right)\right\}.
\label{eqnveqn}
\end{equation}
Expression \eqref{eqnveqn} is monotonic in $\sigma$, and the larger $\sigma$ becomes, the smaller this circular region becomes.  For example, when $\sigma = 2$, the radius of $U$ is approximately $0.5$.  This means the values of $\log\zeta(\sigma + it)$ for $\sigma\geq 2$ are all (approximately) contained in a circular region of radius $.5$ and center at $(0.04,0)$.  For all values of $\sigma > 2$, the absolute value of $\log\zeta(s)$ is smaller than $0.504 < \pi/2$.  Consider $\zeta(s) = - a + i\gamma$ for positive $a$ somewhere in the half-plane $\sigma > 1$.  Then $\log( - a + i\gamma) = \ln| - a + i\gamma| + i(\theta + 2n\pi)$ with $|\theta + 2n\pi| > \pi/2$.  But $\pi/2>1/2$ and so $\text{Re}\,\zeta(\sigma+it)$ is never negative for $\sigma>2$.  In order for the real part of $\zeta(\sigma+it)$ to be negative we must have
\[
\frac {1}{2}\log\frac {\zeta^2(\sigma)}{\zeta(2\sigma)}\geq \pi/2,
\]
or $\sigma< 1.197$.  Once $\sigma$ becomes smaller than this value, the Titchmarsh circles extend beyond the $\pi/2$ line and so $\text{Re}\,\zeta(s)$ can and does become negative.  For the purpose of this study, we only investigate $\zeta(1+it)$ since by continuity, if the function has real part negative on the line $\sigma=1$, then it is also negative for some neighborhood $1<\sigma<1+\epsilon$. 
\section{Initial Considerations}
Consider the expression:
\begin{equation}
-\sum_{n=1}^{N} \log\big(1-p_n^{-\sigma}e^{-i t}\big)
\label{eqn00a}
\end{equation}
for a constant $\sigma$.  What value of $t$ causes the fewest terms of this sum to reach an imaginary part outside the interval $(-\pi/2,\pi/2)$?  It's easy to see that one such value is $t=\pi/2$ since in that case, the terms in \eqref{eqn00a} become $\log\left(1+i \frac{1}{p^{\sigma}}\right)$ and thus give rise to the largest argument. Table \ref{table1} summarizes how this build-up of arguments proceeds.  By the time $N=14$, the imaginary component has exceeded $|\pi/2|$.  Now consider the partial sum: 
$$-\sum_{n=1}^{N}\log\left(1-p_n^{-1} e^{-it\log(p_n)}\right).$$
When $p=2$ and $t=\frac{\pi}{2\log(2)}$ for example, the imaginary part of this term is at its maximum value or $\arctan(1/2)\approx 0.463$.  And in the case of a value of $t$ which simultaneously reduces $t\log(2)$ and $t\log(3)$ $\mod 2$ to a value close to $\pi/2$, the net imaginary part of these two terms would be approximately $0.79$. This is the case for example, when $t=\frac{3186 \pi }{\text{Log}\left[\frac{3}{2}\right]}$.  And when $t$ is larger still, additional "confluences" of this sort can push the first $n$ terms of the sum even closer to a negative imaginary part.  
\begin{table}[h!b!p!]
\caption{Table (default)}
\begin{tabular}{|l|l|}
\hline
N & $\text{Im}\log(1+i p_n)$ \\
\hline
1 & 0.46365 \\
 2 & 0.78540 \\
 3 & 0.98279 \\
 4 & 1.12469 \\
 5 & 1.21535 \\
 6 & 1.29212 \\
 7 & 1.35088 \\
 8 & 1.40346 \\
 9 & 1.44691 \\
 10 & 1.48138 \\
 11 & 1.51363 \\
 12 & 1.54065 \\
 13 & 1.56503 \\
 14 & 1.58829 \\
 15 & 1.60956 \\
\hline
\end{tabular}
\label{table1}
\end{table}  

These considerations do not guarantee the infinite sum will have an imaginary component outside the range $(-\pi/2,\pi/2)$ but does give us some indication how the real part of zeta might grow negative in this part of the half-plane: a value of t which when reduced modulo $2\pi$ for the first $n$ terms, builds up to an imaginary part greater than $\pi/2$ that is simultaneously not offset by the remaining terms in the sum.
The algorithm described in this paper searches for values of $t$ which maximizes this  build-up of initial argument for the first $n$ set of primes, then searches this set of points for values that give rise to a negative real part of $\zeta(1+it)$. 
%
%
\section{Description of the Algorithm}
\label{sectionk}
All algorithms for this study were written in Mathematica 7.0.  The code locates values of t which approximately satisfy \eqref{eqn005} for a set of primes within a tolerance of $\delta$.  One might initially guess that this would reduce to nothing more than a brute-force search for values meeting the tolerances.  However, the modular space of \eqref{eqn005} contains a marvelous symmetry that allows us to locate approximate solutions to $\eqref{eqn005}$ in much less time than a brute-force search.  The solutions to $\eqref{eqn005}$ in this paper are called "confluence points" and denoted by $C_n(\theta,\delta)$ representing a solution to the first $n$ equations of \eqref{eqn005} for particular values of $\theta$ within $\pm\delta$.  In this paper, sometimes when referring to these points, only $C_n$ or "confluence order" is used or if $\theta$ is implied, only $C_n(\delta)$ is stated.
\begin{equation}
M_n=\begin{cases}
t\log(2)&\equiv D\left(\frac{\pi}{2},\delta\right)\mod(2\pi)\\
&\vdots\\
t\log(p_{n})&\equiv D\left(\frac{\pi}{2},\delta\right)\mod(2\pi)
\end{cases}
\label{eqn005}
\end{equation}
A value of $t$ which satisfies the first $n$ equations in \eqref{eqn005} will push the imaginary term of the log sum into a region outside the interval $(-\pi/2,\pi/2)$.  If the net effect of the remaining terms of the sum do not appreciably build-up in the opposite direction, then this value of $t$ will cause $\text{Re}\,\zeta$ to become negative.  The plan then is to locate a large number of $C_n$ points and  check the value of $\zeta(1+it)$ in a neighborhood of each to determine if its real part is negative.  We search along the line $\sigma=1$ since the greater $\sigma$ is, the smaller the argument becomes for each term in the series and the less likely the initial terms will have any affect on the final sum.

\section{Mapping the modular space of $M_n$}

The algorithm described below is designed around a "base confluence set".  We start by determining the confluence points for two primes.  The following description explains a base set in terms of the first two primes.  The algorithm is conceptually the same for at least the pair $(5,7)$ and likely so for other pairs.

\begin{figure}
	\centering
		\includegraphics{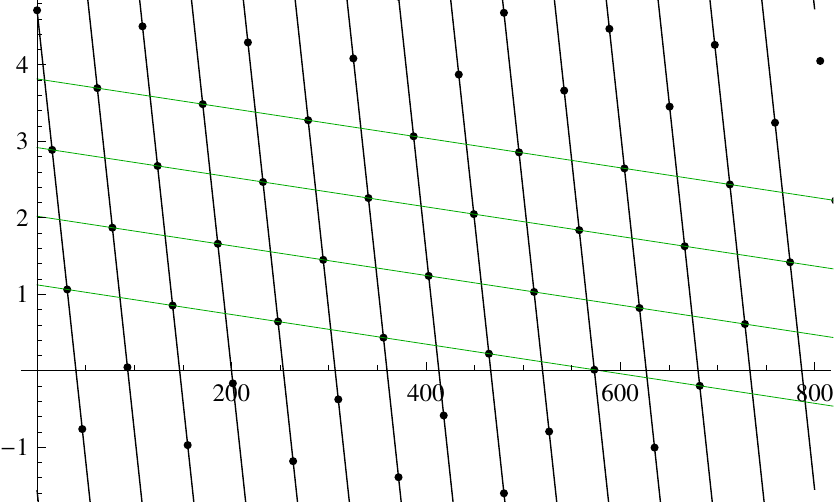}
	\caption{Diagram of $L_n$ and $G_n$}
	\label{fig:modplot2}
\end{figure}

Consider the modular equations
\begin{figure}
	\centering
		\includegraphics{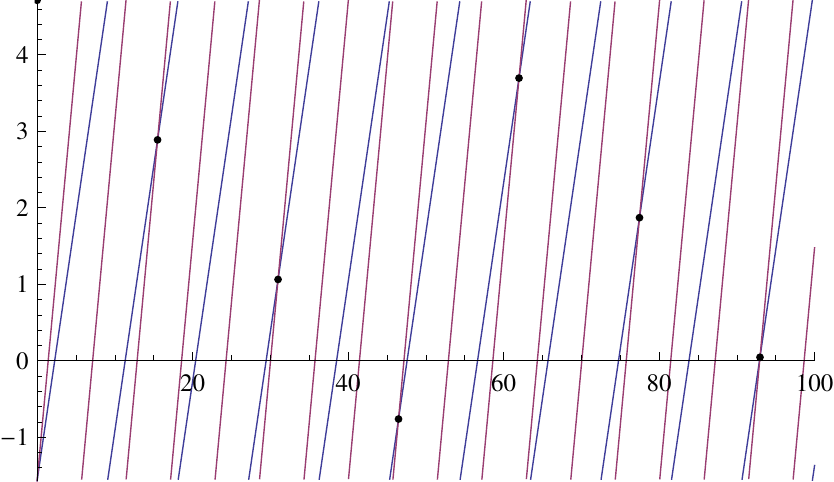}
	\caption{Modulus space $M_n$}
	\label{fig:modplot1}
\end{figure}
%
\begin{equation}
\begin{aligned}
m_2(t)&=\mod(t\log(2),2\pi)-\pi/2\\
m_3(t)&=\mod(t\log(3),2\pi)-\pi/2,
\label{eqn0010}
\end{aligned}
\end{equation}
the plots of which are show in Figure \ref{fig:modplot1} with $m_2(t)$ in blue and $m_3(t)$ in red.  Where the two lines cross, we have a solution to the system:
%
%
\begin{equation}
\begin{aligned}
t\log(2)&\equiv x\mod(2\pi)-\pi/2\\
t\log(3)&\equiv x\mod(2\pi)-\pi/2
\label{eqn0011}
\end{aligned}
\end{equation}
%
%
for some $-\pi/2\leq x<3\pi/2$. One immediately notices a symmetry in Figure \ref{fig:modplot1}:  the confluence points, shown as black dots, line up.  If we connect the dots with diagonal black lines running from the upper left corner to the lower right corner and expand the range of $t$, we obtain the plot in Figure \ref{fig:modplot2} (to avoid clutter, the modular $m_n$ lines have been removed).  Equations for the black lines are given by  
$$L_n(t)=(3\log(2)-2\log(3))t+2n\pi-\pi/2,\quad n=1,2,3,\cdots.$$
The lines $L_n(t)$ serve as a first approximation for locating the next order of confluence points:  we simply calculate where the black lines cross the $t$-axis, check if a base confluence point is near-by, then check for a $C_3$ match within a specified tolerance.  However, the $L_n$ lines are too-closely spaced to offer an efficient algorithm for locating a large number of $C_3$ points quickly.  We can improve the algorithm by considering a second set of symmetries in Figure \ref{fig:modplot2} which is shown as the green lines represented by the equations:
$$G_n(t)=1/3(8\log(2)-5\log(3))t-2n\pi/3+3\pi/2,\quad n=1,2,3,\cdots$$
\begin{figure}
	\centering
		\includegraphics{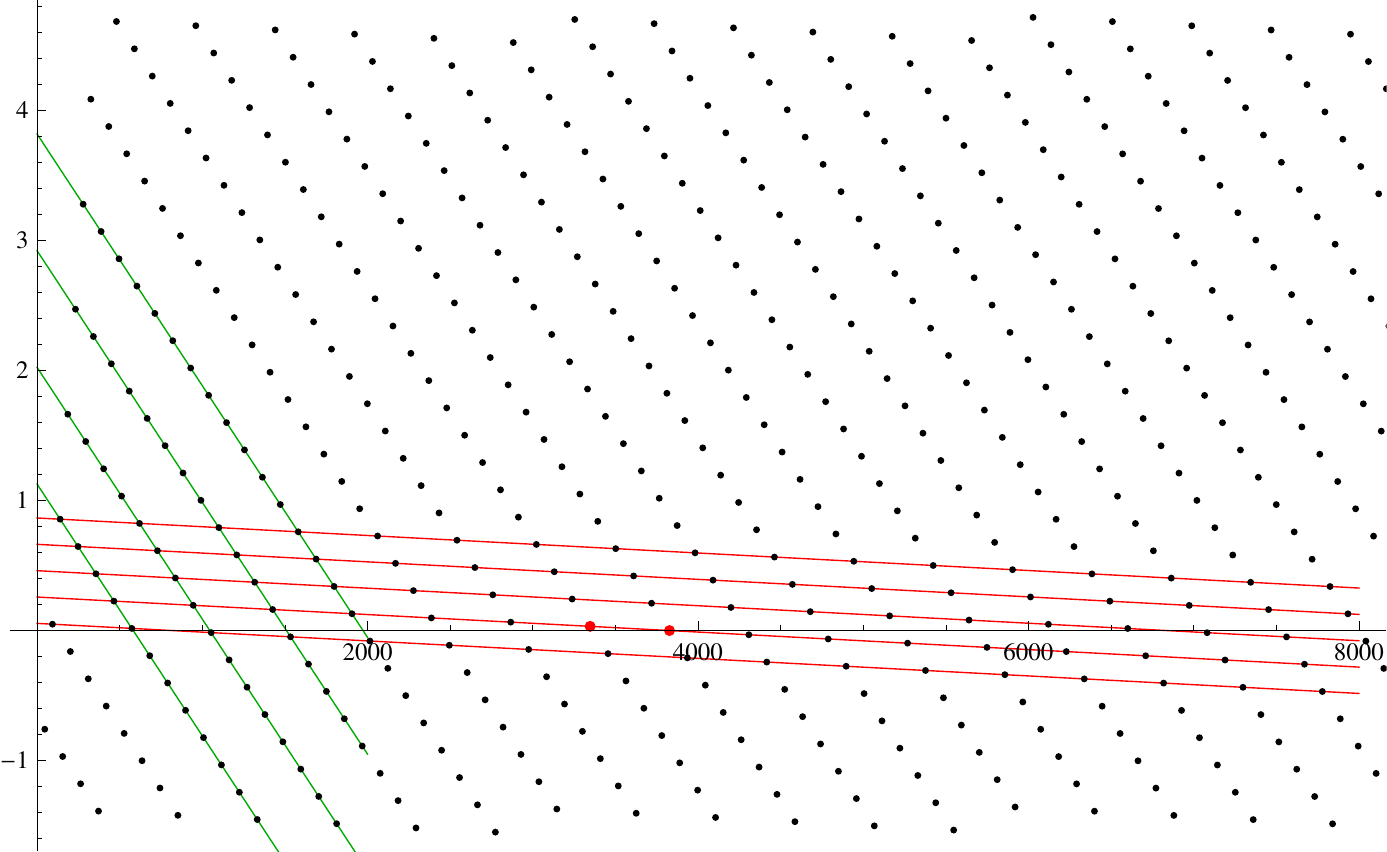}
	\caption{Diagram of $G_n$ and $R(n,m,t)$}
	\label{fig:modplot3}
\end{figure}
In order to further reduce the number of points that need to be checked, we can again form a system of diagonal lines connecting the base confluence points.  These are shown as the red lines in Figure \ref{fig:modplot3} where this time, both the mod lines $m_n$ and black lines $L_n$ were removed for clarity.  The red lines are given by the equations:
$$R(n,m,t)=\log\left(4/3(2/3)^{22/31}\right) t+\pi/62(16 n+4m-47)$$
for $m=1,2,3,4$ and $n=0,1,2,\cdots$.
And although we could continue this process to obtain even further improvements, we stop at the $R(n,m,t)$ lines for this study.  Note in Figure \ref{fig:modplot3},  the $R$ lines cross the $t$-axis with a very small slope with the two nearest $C_2$ confluence points along each red line being very close to the $t$-axis.  Where each $R$ lines crosses the $t$-axis, the nearest-neighbor set denoted by $B_{(2,3)}$, satisfy
\begin{equation}
B_{(2,3)}(\pi/2,\delta)=\begin{cases}
t\log(2)\equiv D\left(\frac{\pi}{2},\delta\right)\mod(2\pi)\\
t\log(3)\equiv D\left(\frac{\pi}{2},\delta\right)\mod(2\pi)
\end{cases}
\end{equation}
for small $\delta$.  An example nearest-neighbor set for one $R$ line is show in the figure as the two red dots.  It is not difficult to determine the zeros of the lines $R(m,n,t)$ and then subsequently,the two nearest neighbors.  This allows for an efficient search of the modular space for the points $C_2(\pi/2,\delta)$ requiring only a few computations for every $2000$ increment in $t$.
Each confluence point $C_2$, along the $R$ lines is at a distance of $\Delta_z(2,3)=\frac{62\pi}{\log(3)-\log(2)}\approx 480$.  The distance between each successive zero of the $R$ lines is $\frac{2\pi}{53 \log(3)-84 \log(2)}\approx 3009$ and the slope of each $R$ line is $m_R=\log(4/3(2/3)^{22/31})\approx -6.7\times 10^{-5}.$  Thus, knowing the zero point of each $R(m,n)$ and $\Delta_z(2,3)$, we can determine with relatively few calculations, the nearest $C_2$ confluence points closest to the t-axis and then check only these for higher-ordered confluence points.  But it turns out that we can do much better than this because the pattern of confluence points in one range of values of $t$ goes a long way in determining the pattern in a following range of values.  Some drift is of course expected but in actual runs of the software, this turned out to be surprising little over a very large interval. We can exploit this property by using a set of lower-ordered confluence points as a map for the next higher set of confluence points thereby greatly reducing the time to locate these points.

It's easy to see that the maximum height reached by either the previous or next $C_2$ point along each $R_n$ line is $m_R \frac{62\pi}{\log(3)-\log(2)}\approx 0.032$.  Therefore, if the tolerance $\delta$ is set greater than $0.032$, than both the previous and next point will be below this value and can be chosen without checking the residue mod $2\pi$.

The following Mathematica code implements this approach to locate the base-confluence points $B_{(2,3)}$.
$$
\begin{aligned}
&\hspace{-15pt}\texttt{For[n=1,n$\leq$nMax,n++,} \\
&\hspace{-10pt}\texttt{For[m=1,m$\leq$4,m++,} \\
&\texttt{zeropt}=\frac{(-31+4m+16n)\pi}{2(53\log(3)-84\log(2))}\texttt{;}\\
&\texttt{startpt}=\frac{2(31-7m+3n)\pi}{\log(3)-\log(2)}; \\
&\texttt{intpt}=\text{IntegerPart}[\frac{62}{k}\left(\texttt{zeropt}-\texttt{startpt}\right)];\\
&\texttt{prevpt}=\left(\texttt{startpt}+\frac{62\pi\texttt{inpt}}{k}\right);\\
&\texttt{nextpt}=\left(\texttt{prevpt}+\frac{62\pi}{k}\right);\\
&\texttt{If[theDegree$<$maxDifference,} \\
&\hspace{15pt}\texttt{thePointList23=Append[thePointList23,prevpt];} \\
&\hspace{15pt}\texttt{thePointList23=Append[thePointList23,nextpt];} \\
&\texttt{,}\\
&\hspace{15pt}\texttt{testprev=Mod[prevpt Log[2],2$\pi$]-$\pi/2$;}\\
&\hspace{15pt}\texttt{testnext=Mod[nextpt Log[2],2$\pi$]-$\pi/2$;}\\
&\hspace{20pt}\texttt{If[0 $<$ testprev $<$ theDegree,}\\
&\hspace{25pt}\texttt{thePointList2=Append[thePointList2,prevpt];}\\
&\hspace{20pt}\texttt{];}\\
&\hspace{20pt}\texttt{If[0 $<$ testnext $<$ theDegree,}\\
&\hspace{25pt}\texttt{thePointList2=Append[thePointList2,nextpt];}\\
&\hspace{20pt}\texttt{];}\\
&\texttt{];}\\
&\hspace{-10pt}\texttt{];}\\
&\hspace{-15pt}\texttt{];}
\end{aligned}
$$
From this starting set of points, we go on to find the set of $C_3$ points by calculating the gaps between successive $C_2$ points and then use those gaps to predict with good success, the location of the next $C_3$.  In particular, if we use the above code to locate the first $10,000$ $C_2(\pi/2,0.01)$ points, we find only three possible gap sizes:
$$\left\{g_1,g_2,g_3\right\}=\left\{\frac{778\pi}{\log(3/2)},\frac{7360\pi}{\log(3/2)},\frac{8138\pi}{\log(3/2)}\right\}$$We can then use this table to search for $C_3$ point by simply checking three possible locations.  For example, if a $C_2$ point was located at $t_2$, then we need only check the three points $t_2+g_1$, $t_2+g_2$, $t_2+g3$. 
The following code implements the next stage of the search by locating the $C_3$ points using the previously-calculated $C_2$ gap table.  In the code, we have a While loop which checks for the $C_2$ points against the corresponding gap table.  As the confluences grow in order, we add conditionals for additional primes.  The following code is only the main routine.  In actual practice, pre-processing of the previous confluence data is necessary.  The variable \texttt{mintval} is then initially set to the last $C_{n}$ point (if any) found in the pre-processed data. After each acquisition, the data must then be post-processed into a sorted gap table for the next confluence search.
$$
\begin{aligned}
&\texttt{thePointList3}=\{\};\\
&\texttt{mintval=lastin23};\\
&\hspace{10pt}\texttt{n=1}; \\
&\texttt{For mypoints=1},\texttt{mypoints}\leq \texttt{maxN,mypoints++} \\
&\hspace{10pt}\texttt{n=1}\\
&\hspace{10pt}\texttt{testpt}=\texttt{mintval}+\texttt{gaptable23}[[n]];\\
&\hspace{10pt}\texttt{While[} \\
&\hspace{15pt}\texttt{Not[0$\leq$ Mod[test Log[2],$2\pi$]-$\pi/2<$theDegree]},\\
&\hspace{20pt}\texttt{n++};\\
&\hspace{20pt}\texttt{If[n>}\texttt{gaptable23len},\\
&\hspace{25pt}\texttt{Print["Gap table 23 exhausted"]};\\
&\hspace{25pt}\texttt{Abort[];}\\
&\hspace{25pt}\texttt{,}\\
&\hspace{25pt}\texttt{testpt}=\texttt{mintval}+\texttt{gaptable23[[n]]};\\
&\hspace{25pt}\texttt{];}\\
&\hspace{20pt}\texttt{];}\\
&\hspace{15pt}\texttt{mintval}=\texttt{testpt};\\
&\hspace{15pt}\texttt{If[Mod[mintval Log[5],2$\pi$]-$\pi/2$<theDegree,}\\
&\hspace{15pt}\texttt{thePointList3}=\texttt{Append[thePointList3,mintval]};\\
&\hspace{10pt}\texttt{];}\\
&\hspace{5pt}\texttt{];} \\
\end{aligned}
$$
With this code and some additional pre-processing and post-processing code, we are able to find a large number of $C_3(\pi/2,\delta)$ points quickly.  For example, we can check about  50,000 $C_2$ points in about $40$ seconds with the machine this code was tested on.  We next calculate the gap table for a sufficient number of $C_3$ points.  In actual practice, this was usually around $1500$ points.  In one run of the $C_3$ code, we obtained the following gap table:

$$\texttt{gaptable3}=\left\{\frac{91260\pi}{\log(3/2)},\frac{360014\pi}{\log(3/2)},\frac{451274\pi}{\log(3/2)},\frac{921936\pi}{\log(3/2)},\frac{1013196\pi}{\log(3/2)},\frac{1373210\pi}{\log(3/2)},\frac{1464470\pi}{\log(3/2)}\right\}$$ 
Notice the gap table has grown in size from the first table.  This is quite expected considering we are working with residues $\mod 2\pi$.  The next step of course is to use $\texttt{gaptable3}$ to locate the $C_4(\pi/2,\delta)$ points.  We then add another conditional to the While loop to check for a $C_3$ confluence point.  We then create a gap table of those points and continue the search and if necessary, adjust the tolerance.  The routine is coded to terminate a run if after checking all values suggested by the gap table for a particular point, the routine fails to find the desired confluence point at the specified tolerance.  In this case, a small change to the tolerance was found to correct for the drift at least in the range studied.
\section{Acquisition of Data}
All code was run on a $2.2$ GHz machine with one Gb RAM. In order to study this algorithm, three runs were performed.  The first two were run with a base-confluence set $B_{(2,3)}$.  The details of constructing this set were described in the previous sections.  The third set was done with a base-confluence set $B_{(5,7)}$, and it's construction is virtually identical to the $B_{(2,3)}$ set except the equations for $L_n,G_n,$ and $R_n$ are different reflecting the use of different primes in their construction. In all cases, $\theta=\frac{\pi}{2}$.
%
%
\subsection{Run 1}

Using base-confluence $B_{(2,3)}$, and $\delta$ set to 0.1, the algorithm located approximately $1000$ points for each confluence up to $C_{10}$ in approximately $90$ minutes. The variable \texttt{nMax} was set to $50,000$. This means that $50,000$ $C_{n-1}$ terms were checked in order to obtain the $C_n$ set within the currently set tolerance.  In all cases, the tolerance was "globally" set meaning that each modular calculation was checked against the same value of $\delta$. Table \ref{table1} summarizes the result of this run.  "First" and "Last" are the first and last points of the set, "Gap" is the average distance between $C_n$ points as determined by this algorithm, and $k=\frac{\pi}{\log(3)-\log(2)}$, that is, the first $C_3$ point is located at $\frac{8274\pi}{\log(3)-\log(2)}$.

\begin{table}[ht]
\caption{Run 1}
\centering
\begin{tabular}{|c|c|r|r|c|}
\hline
Type & Points & First & Last & Gap \\
\hline
$B_{(2,3)}(.1)$ & 7992 &  $432k$ & $1551976k$ & $10^3$ \\
$C_{3}(.1)$ & 923 &   $8274k$ & $9342654k$ & $10^5$ \\
$C_{4}(.1)$ & 817 &   $171406k$ & $450025676k$ & $10^6$ \\
$C_{5}(.2)$ & 1601 &  $6538724k$ & $11407973800k$ & $ 10^8$ \\
$C_{6}(.2)$ & 1625 & $180883380k$  & $101681208846k$ & $ 10^9$ \\
$C_{7}(.2)$ & 1659 & $1524427526k$ & $556034603992k$ & $10^9$ \\
$C_{8}(.2)$ & 1689 & $34995876276k$ & $3927590605616k$ & $10^{10}$ \\
$C_{9}(.3)$ & 2438 & $232692763660k$ & $20585726758472k$ & $10^{11}$ \\
$C_{10}(.4)$ & 3197 & $5206151247198k$ & $53443637744022k$& $ 10^{11}$ \\
\hline
\end{tabular}
\label{table1b}
\end{table}  
Note the tolerance was relaxed at $C_5$ and $C_9$ due to gap table exhaustion. Figure \ref{fig:run1c6} represents a "confluence portrait" for one $C_6$ value obtained in this run and serves to graphically illustrate the solutions for $M_6$ of \eqref{eqn005}.  Each line in the figure represents one equation of \eqref{eqn005} plotted in the range $(t-\alpha,t+\alpha)$ for some $\alpha$.  Notice how all the equations confluence at this $C_6$ point hence the name given to these points.  

Consider the six terms of the partial sum of \eqref{eqn00a} which for simplicity we take $\sigma=1$ since we are only interested in an approximation:

\begin{figure}
	\centering
		\includegraphics{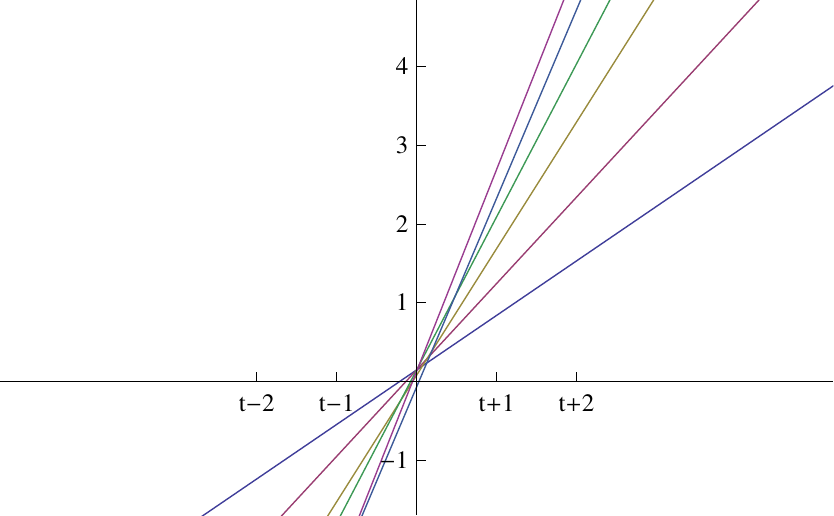}
	\caption{$C_6(\pi/2,0.2)$ at $t_6=20\,430\,730\,768 k$ }
	\label{fig:run1c6}
\end{figure}

$$-\sum_{n=1}^{6}\log\left(1-\frac{1}{p} e^{-i t_6 \log(p)}\right)\approx -0.126968 - 1.32215 i, $$
and note the imaginary part is already close to $|\pi/2|$ and even although this gives no indication of how the remaining terms will contribute to the sum, it illustrate how the confluences are at least affecting the first few terms of the log sum.  
%
%
\subsection{Run 2}

\begin{figure}
	\centering
		\includegraphics{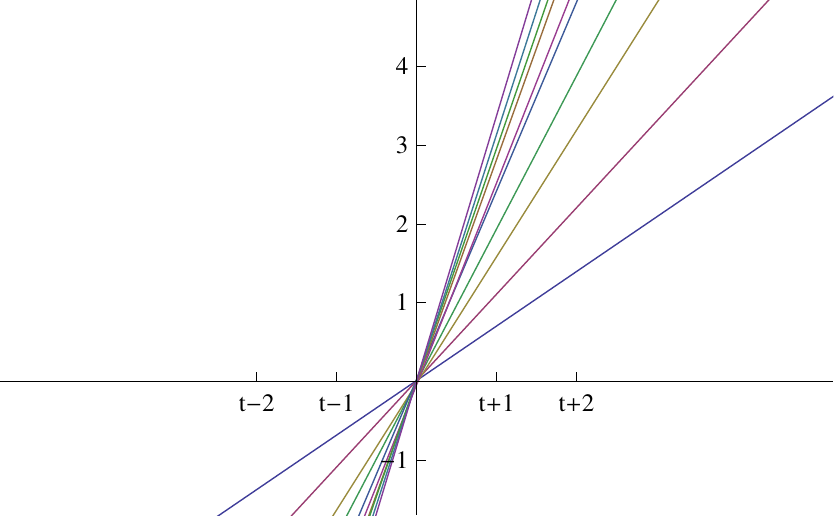}
	\caption{$C_{10}(\pi/2,0.08)$ at $t=34\,494\,360\,555\,864\,694k$}
	\label{fig:c10run2}
\end{figure}

The second set of data was initially acquired at $\delta=0.05$ (except the base-set) and with a sampling size of $200,000$.  The gap tables successfully located on average $1600$ confluence points until $C_9$ was reached.  Then $\delta$ was relaxed to $0.07$ and the $C_9$ set successfully acquired.  Alternatively, we could have re-acquired a larger gap table for $C_8$.  The acquisition of $C_{10}$ exhausted the $C_9$ gap table.  This set of confluence points was re-acquired at $\delta=0.08$ until after the acquisition of $1273$ terms, the gap table was again exhausted.  Total execution time was approximately five hours.  Figure \ref{fig:c10run2} shows a typical $C_{10}$ confluence portrait at this value of $\delta$.  Note the equations at the origin confluence at a smaller tolerance compared to those in Figure \ref{fig:run1c6}.

\begin{table}[ht]
\caption{Run 2}
\centering
\begin{tabular}{|c|c|r|r|c|}
\hline
Type & Points & First & Last & Gap \\
\hline
$B_{(2,3)}(.01)$ & 1239 & $3186k$ & $1551556k$ & $10^4$ \\
$C_{3}(.05)$ & 1599 &    $118922k$ & $242220222k$ & $10^6$ \\
$C_{4}(.05)$ & 1602 &   $66656806k$ & $19751098706k$ & $10^8$ \\
$C_{5}(.05)$ & 1642 &  $1095953254k$ & $1684031570650k$ & $10^{10}$ \\
$C_{6}(.05)$ & 1624 &   $10971337166k$  & $46958108569642k$ & $10^{11}$ \\
$C_{7}(.05)$ & 1563 &  $5837564637802k$ & $963484780912880k$ & $10^{12}$ \\
$C_{8}(.05)$ & 1608 &  $235403791542126k$ & $41141345684941152k$ & $10^{14}$ \\
$C_{9}(.07)$ & 2199 &  $1409137033297936k$ & $838179310625884816k$ & $10^{15}$ \\
$C_{10}(.08)$ & 1273 & $3449430555864694k$ & $5206649605555820266k$& $10^{16}$ \\
\hline
\end{tabular}
\label{table2}
\end{table}  

%
%
\subsection{Run 3}

\begin{figure}
	\centering
		\includegraphics{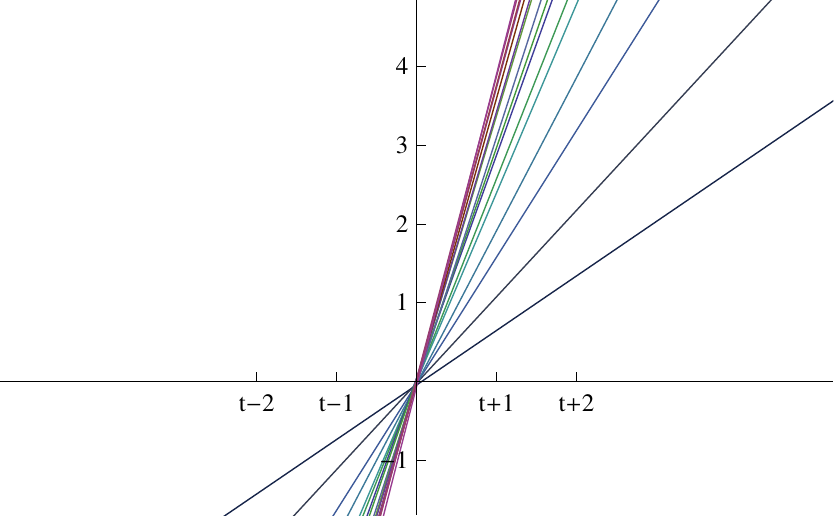}
	\caption{$C_{16}(\pi/2,0.063)$ at $t=12253527959225463513391519458k$}
	\label{fig:firstc16}
\end{figure}

Since one objective of this study was to devise an efficient algorithm for numerically solving the system of congruences defined in \eqref{eqn005}, one would like to know how well it worked for extremely large values of $t$. This third run attempted to determine this by running the algorithm at an initial tolerance of $\delta=0.01$ and sampling size of $500,000$.  As Table \ref{tabler3} suggests, the algorithm was stable enough to acquire $1765$ $C_{14}$ points with an average gap interval of $10^{26}$ at a tolerance of $0.05$.  This allowed the acquisition of $121$ $C_{15}$ before exhausting the $C_{14}$ gap table which by then had grown to about $2000$ entries.  This data took approximately 24 hours to acquire.  The final $121$ $C_{15}$ points were checked for a minimum $C_{16}$ confluence which was found at a tolerance of $0.063$ and is shown in Figure \ref{fig:firstc16}.  In Table \ref{tabler3}, $k=\frac{\pi}{\log(7)-\log(5)}$.
\begin{table}[ht]
\caption{Run 3}
\centering
\begin{tabular}{|c|c|r|r|c|}
\hline
Type & Points & First & Last & Gap \\
\hline
$B_{(2,3)}(.01)$ & 835 & $1188k$ & $105918k $ & $10^3$ \\
$C_{3}(.01)$ & 239 & $415346k$ & $188449608k $ & $10^5$ \\
$C_{4}(.01)$ & 1760 & $951398k$ & $341342656k $ & $10^6$ \\
$C_{5}(.01)$ & 785 & $945400159026k$ & $248625146191138k $ & $10^{12}$  \\
$C_{6}(.01)$ &  782 & $358872295235106k$ & $53873791685357248k $ & $10^{15}$ \\
$C_{7}(.01)$ &  2449 & $25465934548192130k$ & $24901118319251218470k $ & $10^{17}$\\
$C_{8}(.01)$ &  1607 & $4228669878993511528k$ & $3026698278169709217444k $ & $10^{19}$ \\
$C_{9}(.01)$ &  1581 & $443049036634051704488k$ & $106247940250479965840452k $ & $10^{20}$\\
$C_{10}(.02)$ & 1608 &  $15784886571659169335556k$ & $2201816585968897375148202k $ & $10^{22}$ \\
$C_{11}(.03)$ & 2418 & $325647381902527914552536k$ & $37724324922892546494833526k$  & $10^{23}$\\
$C_{12}(.03)$ & 1800 & $4833028477320576769793928k$ & $502072051726121146180272476k$ & $10^{24}$ \\
$C_{13}(.04)$ & 3924 & $60624317376187697337278868k$ & $5796475948831968271200664166k$ & $10^{25}$ \\
$C_{14}(.04)$ & 1765 & $510082174019371223753727784k$ & $19898784954288834658068768846k$ & $10^{26}$ \\
$C_{15}(.05)$ & 121 & $6450267346766950732885822366k$ & $26894796143731914838108040998k$ & $10^{27}$ \\
\hline
\end{tabular}
\label{tabler3}
\end{table} 

\section{Processing the Data Sets}

For the first run, we tested each entry in the confluence tables up to $C_6$ by calculating $\zeta(1+i(t\pm1))$ at $20$ equally-spaced points in the interval and obtained a total of ten regions along the line $\text{Re}(s)=1$ for which the real part of zeta became negative.  Sample points in the indicated $C_n$ region are given in Table \ref{table15}.  The value of zeta in the regions was determined using the Mathematica "N" command using  arbitrary precision arithmetic with the final answer reported to thirty digits of accuracy to the right of the decimal place using the following command for each of the ten points in Table \ref{table15}:
$$\texttt{N[zeta[1+it],$\{\infty,30\}$]}$$
The smallest value of $t$ detected in the data sets described above in which $\text{Re}\,\zeta(1+it)<0$ was found in the vicinity of $3.4\times 10^7$ as shown in Table \ref{table15} (recall the factor $k$).  All the values in Table \ref{table15} had negative imaginary part due to the search criteria selecting $\pi/2$ residues.  No points in the first four confluence tables of the remaining runs were found to give rise to a negative real part of $\zeta$ at least using the method described above.  The remaining confluence tables were not inspected at this time because of the impracticality of actually computing the zeta function (on the machine used to analyze the data) at the relatively large values of $t$ in these tables.
\begin{table}[h!b!p!]
\caption{Points where $\text{Re}\,\zeta(1+it)<0$}
\begin{tabular}{|c|c|l|l|}
\hline
ID & Type &t & $\zeta(1+it)$ \\
\hline
1 &$C_3$ & $4\,378\,640k\phantom{111k}-2/5$ & $-0.009-1.22i$\\
2 & $C_4$ & $415\,782\,314k\phantom{1k}-2/5$ & $-0.024-1.23i$\\
3 & $C_6$ & $20\,430\,730\,768k-1/10$ &$-0.015-1.08i$ \\
4 &$C_6$ & $25\,705\,015\,862k-1/5$ &$-0.002-1.20i$ \\
5 & $C_6$ & $47\,668\,373\,108k-1/10$ & $-0.027-0.96i$ \\
6 & $C_6$ & $53\,761\,507\,682k$ & $-0.009-1.00i$ \\
7 & $C_6$ & $62\,484\,882\,686k-1/5$ & $-0.013-1.13i$ \\
8 & $C_6$ & $65\,421\,460\,042k-1/10$ & $-0.003-0.91i$ \\
9 & $C_6$ & $97\,190\,286\,104k-1/10$ & $-0.018-1.00i$ \\
10 & $C_6$ & $99\,154\,858\,182k-3/10$ & $-0.026-1.25i$\\
\hline
\end{tabular}
\label{table15}
\end{table}  

\section{Conclusions}

One can easily locate a value of $t$ where $\text{Re}\,\zeta(1+it)<0$ by trial an error by simply  sampling the $t$-axis at intervals of $1/10$ starting at $t=0$.  Doing so, one first finds a negative real part at $t=682112.9$. However, this study attempted to predict where the real part of the function would turn negative without a trial and error approach.  A surprising result to come from this study is the relative stability of the method.  One would have guessed that such a method would quickly degrade due to chaotic drift in the modular space rendering it virtually impossible to predict any future values let alone values up to $10^{29}$.  It would be interesting in a further study to try and explain what is the source of this stability.  

Additionally:
\begin{enumerate}
\item
What are the precise shapes of these negative contours where $\text{Re}\,\zeta(\sigma+it)<0$ and how do they change as we go up the $t$-axis?  
\item
This study focused on two base-confluence algorithms:  $(2,3)$ for the first two trials, and $(5,7)$ for the third trial.  It is not known how the algorithms are affected by different base prime pairs.

\item
This study demonstrated some structure in the modular space of \eqref{eqn005}.  Is there a similar periodicity in the domains where $\text{Re}\,\zeta(\sigma+it)$ dips below the $\sigma-t$ plane?  Can one predict where such domains might be located analogous to how the various confluence points were "predicted" using the gap tables? 

\item 
This study did not attempt to investigate in detail the drift of values in the modular space.  When a gap table became exhausted, the tolerance was simply increased or a greater number of $C_{n-1}$ values were located and a new larger gap table created.  It would be interesting to better understand this drift.
\end{enumerate}

\vspace{20pt}
\address{\noindent D.C. Milioto\\LaPlace, LA. 70068 (U.S.A)  \\ email: miliotodc@rtconline.com}


\begin{thebibliography}{9}                                                                    
\bibitem {Edwards}\textsc{H.Edwards}, \textit{Riemann's Zeta Function}, Academic Press, New York,1974.
\bibitem {Flanigan}\textsc{F.J.Flanigan},\textit{Complex Variables, Harmonic and Analytic Functions}, Dover Publications, Inc., New York, 1972.

\bibitem {Grenne}\textsc{R.E. Greene}, \textit{Function Theory of One Complex Variable}, American Mathematical Society, Providence, R.I., 2006.

\bibitem {Ingham} \textsc{A.E. Ingham},\ \textit{The Distribution of Prime Numbers}, Cambridge University Press, Cambridge, Mass., 1995.

\bibitem{Marsden} \textsc{J.F. Marsden and M.J. Hoffman},\ \textit{Basic Complex Analysis},W. H. Freeman, New York, 1987.

\bibitem{Stein} \textsc{E.M. Stein and R. Shakarchi},\textit{Complex Analysis}, Princeton University Press,Princeton, N.J., 2003.

\bibitem {Titchmarsh2} \textsc{E.C. Titchmarsh},\ \textit{The Theory of the Riemann Zeta Function}, Cambridge University Press,New York, 1987.

\end{thebibliography}
\end{document}